\documentclass{amsart}
\usepackage[nosections]{macs-ddc}
\usepackage{amsrefs}
\usepackage[all]{xy}
\CompileMatrices
\SelectTips{cm}{}
\allowdisplaybreaks[2]

\date{12 August 2002}

\begin{document}

\title{The Dixmier-Douady Class of Groupoid Crossed Products}

\author{Paul S. Muhly}
\address{Department of Mathematics\\
University of Iowa\\
Iowa City, IA 52242}
\email{muhly@math.uiowa.edu}

\author{Dana P. Williams}
\address{Department of Mathematics\\
Dartmouth College\\
Hanover, NH 03755-3551}
\email{dana.williams@dartmouth.edu}

\begin{abstract}
We give a formula for the Dixmier-Douady class of a continuous-trace groupoid
crossed product that arises from an action of a locally trivial, proper,
principal groupoid on a bundle of elementary $\cs$-algebras that
satisfies Fell's condition.
\end{abstract}

\maketitle

\section{Introduction}

Throughout this note $G$ will denote a second countable, locally compact,
principal groupoid with Haar system,  $\A$ will denote an elementary
$C^{*}$-bundle over $\go$ that satisfies Fell's condition, and $\alpha$
will denote a continuous action of $G$ on $\A$ via isomorphisms.
Thus the pair $(\A,\alpha)$ is exactly
what is needed to define an element in the Brauer group $\brg$ as
defined in \cite{kmrw:ajm98}*{Definition~2.14}.  As a special case of
\cite{fmw:xx}*{Theorem~1}, it follows that the groupoid crossed
product $\cs(G;\A)$ has continuous-trace if and only if $G$ is a
proper principal groupoid.  Thus if $\cs(G;\A)$ has continuous trace
we can, and do, assume that that $G=\Omega\pstar\Omega = \set{(\omega,
  \omega')\in\Omega\times\Omega:p(\omega)=p(\omega')}$ for a
continuous open surjection $p:\Omega\to Y$.  In this note, we want to
consider its Dixmier-Douady class $\delta\bigl(\cs(G;\A)\bigr)$, and
we compute $\delta\bigl(\cs(G;\A)\bigr)$ when $p$ is locally trivial.
In this case we say that $G$ is a \emph{locally trivial, proper
  principal groupoid}.  Our approach is motivated in part by
\cite{raeros:tams88}*{\S1} where Raeburn and Rosenberg consider the
case where $G$ is the transformation group groupoid $G=H\times\Omega$
with $\Omega$ a locally trivial principal $H$-bundle, and $\alpha$ is
pulled back from a locally unitary action $\gamma$ of $H$ on a stable
continuous-trace \cs-algebra $B$.  Just as in
\cite{raeros:tams88}*{Theorem~1.1}, we show that we can assume that
$\A$ is the pull-back $p^{*}\B$ by the orbit map $p:G\to\gomg$ of some
locally trivial $\K$-bundle\footnote{Throughout, $K$ will denote the
  space of compact operators on a separable, infinite dimensional
  Hilbert space $H$. By a locally trivial $K$-bundle over a space $X$,
  we shall mean a locally trivial fibre bundle over $X$, with fibre
  $K$ and structure group $\Aut(K)$.  We shall use \cite{rw:morita} as
  our basic reference on fibre bundles, sheaves, cohomology, etc.}
over $\gomg$.  In the event $\alpha$ satisfies an additional
hypothesis --- similar to being pulled back from a locally unitary
action as in \cite{raeros:tams88} --- we show that 
$\delta\bigl(\cs(G;\A)\bigr)$ is a naturally defined perturbation of
$\delta(B)$ in complete analogy with
\cite{raeros:tams88}*{Theorem~1.5}.

\section{Locally Trivial Proper Principal Groupoids}

Suppose that $\Omega$ and $Y$ are second countable
locally compact Hausdorff spaces, and that $p:\Omega\to Y$ is a
continuous, open surjection.  Let $G$ be the proper principal groupoid
$\Omega\pstar\Omega$ and identify $\go$ with $\Omega$ and $\gomg$ with
$Y$. We'll say that $G$ is
\emph{locally trivial} if $p:\Omega\to Y$ is a locally trivial fibre
bundle with fibre $X$.  That is, we assume there is a cover
$\U=\set{U_{i}}$ of $Y$ and continuous maps $h_{i}:p^{-1}(U_{i}) \to
U_{i}\times X$ such that the diagram 
\begin{equation*}
  \xymatrix{p^{-1}(U_{i}) \ar[rr]^-{h_{i}} \ar[dr]_-{p}&&U_{i}\times X
  \ar[dl]^{\operatorname{pr}_{1}}\\&U_{i} }
\end{equation*}
commutes for each $i$.  In particular, we assume that
there are continuous functions
$\xi_{i}:p^{-1}(U_{i})\to X$ such that 
\begin{equation}
  \label{eq:4}
  h_{i}(w)=\bigl(p(\omega),\xi_{i}(\omega)\bigr)\quad\text{for all
  $\omega\in p^{-1}(U_{i})$.}
\end{equation}
Consequently, if
$U_{ij}:=U_{i}\cap U_{j}\not=\emptyset$,
then the diagram
\begin{equation*}
  \xymatrix{U_{ij}\times X \ar[dr]^-{h_{i}^{-1}} \ar@{..>}[dd]
    \ar@/^1pc/[drr]^{\operatorname{pr}_{1}}  \\
&p^{-1}(U_{ij}) \ar[r]^-{p} \ar[dl]^-{h_{j}} & U_{ij} \\
U_{ij}\times X \ar@/_1pc/[urr]_{\operatorname{pr}_{1}}}
\end{equation*}
commutes.  This, in turn, implies that
 for each $u\in U_{ij}$ there must be a continuous map
$\phi_{ij}(u):X\to X$ such that 
\begin{equation}
  \label{eq:1}
  h_{j}\circ
  h_{i}^{-1}(u,x)=\bigl(u,\phi_{ij}(u)(x)\bigr)\quad\text{for all
  $(u,x)\in U_{ij}\times X$.}
\end{equation}
Further, since $\phi_{ij}(u)^{-1}=\phi_{ji}(u)$, we see that each
$\phi_{ij}(u)$ is a homeomorphism of $X$. Consequently, $\phi_{ij}$ may
be viewed as a function from $U_{ij}$ to $\Homeo(X)$.
A straightforward computation shows
that
\begin{equation}
  \label{eq:2}
  \phi_{ik}(u)=\phi_{jk}(u)\circ \phi_{ij}(u)\quad\text{for all $u\in
  U_{ijk}:=U_{i}\cap U_{j}\cap U_{k}$.} 
\end{equation}

\begin{remark}
  \label{rem-homeox}
  $\Homeo(X)$ can be made into a topological group in such a way that
  $h_{n}\to h$ in $\Homeo(X)$ if and only if given a net $x_{n}\to x$
  in $X$ (with the \emph{same} index set), then $h_{n}(x_{n})\to h(x)$
  and $h_{n}^{-1}(x_{n})\to h^{-1}(x)$.  Then it not hard to see that
  $\phi_{ij}$ is continuous and that the transition functions
  $\set{\phi_{ij}}$ determine $p$ in the usual way.
\end{remark}

\begin{example}
  \label{ex-principal-bun}
  Of course, the basic example of a locally trivial proper principal
  groupoid is the transformation groupoid associated to a principal (left)
  $H$-space $\Omega$ for a locally compact group $H$.  In this case,
  the fibre space $X$ is just $H$, and we also want the local trivializations
  $h_{i}$ to be $H$-equivariant.  In particular, there are continuous
  functions $s_{ij}:U_{ij}\to H$ such that
  $\phi_{ij}(u)(t)=ts_{ij}(u)$ for all $t\in H$.  Furthermore,
  equation~\eqref{eq:2} is equivalent to
  \begin{equation}
    \label{eq:10}
    s_{ij}(u)s_{jk}(u)=s_{ik}(u)\quad\text{for all $u\in U_{ijk}$.}
  \end{equation}
  Therefore the elements $\set{s_{ij}}$ determine the class $[p]$ of
  the principal bundle $p$ in $H^{1}(Y,\sheaffont{H})$ where
  $\sheaffont{H}$ is the sheaf of continuous $H$-valued functions and
  $H^{1}(Y,\sheaffont{H})$ is the the first sheaf cohomology \emph{set}
  determined by $\sheaffont{H}$ \cite{rw:morita}*{Remark~4.54}.  When we
  return to this example in the sequel, we will identify the
  transformation group groupoid $H\times\Omega$ with
  $G=\Omega\pstar\Omega$ via the map $(t,\omega)\mapsto
  (\omega,t^{-1}\cdot \omega)$.\footnote{Recall that $H\times\Omega$
    is the groupoid with unit space $\set e\times \Omega$ identified
    with $\Omega$ and with range and source maps
    $s(t,\omega)=t^{-1}\cdot \omega$ and $r(t,\omega)=\omega$.  Then
    we have $(t,\omega)(s,t^{-1}\cdot \omega)= (ts,\omega)$ and
    $(t,\omega)^{-1} = (t^{-1},t^{-1}\cdot \omega)$.}
\end{example}

\begin{remark}
  \label{rem-choices}
  It is a matter of taste as to whether $\phi_{ij}$ or $\phi_{ji}$
  appears in \eqref{eq:1}.  Our taste might seem off in view
  of \eqref{eq:2}, but we have purposely endured bitter herbs in order
  enjoy \eqref{eq:10} even if $H$ is not abelian.  In either
  case, it is
  important fact of life that any formula for the Dixmier-Douady class
  that depends on standard topological data, such as transition
  functions like the $\phi_{ij}$, depends up to a
  sign on the choices such as that made in \eqref{eq:1}.  This will be
  important in comparing our result to other calculations in the
  literature 
  (cf. Example~\ref{ex-last}). 
\end{remark}

Now we want to see how the locally triviality of $p$ is reflected in the
groupoid structure of $G$.
We can define a topological groupoid isomorphism
$k_{i}:G\restr{p^{-1}(U_{i})}\to X\times U_{i}\times 
X$ by 
\begin{equation}
  \label{eq:5}
  k_{i}(\omega,\omega'):=
  \bigl(\xi_{i}(\omega),p(\omega),\xi_{i}(\omega')\bigr)  .  
\end{equation}
Here $X\times U_{i}\times X$ is the groupoid that has unit space identified with
$X\times U_{i}$, orbit space identified with $U_{i}$ and
multiplication given by $(x,u,y)(y,u,z)=(x,u,z)$.  
If $\ptwo:G\to
\go/G\cong Y$ is given by $\ptwo(\omega,\omega'):=p(\omega)$, then we   
have commutative diagrams
\begin{equation*}
  \xymatrix{G\protect\restr{p^{-1}(U_{i})} 
\ar[rr]^-{k_{i}} \ar[dr]_-{\protect\ptwo}&&X\times U_{i}\times X
  \ar[dl]^{\operatorname{pr}_{2}}\\&U_{i} }
\end{equation*}
and
\begin{equation*}
  \xymatrix{X\times U_{ij}\times X \ar[dr]^-{k_{i}^{-1}} \ar@{..>}[dd]
    \ar@/^1.5pc/[drr]^{\operatorname{pr}_{2}}  \\
&G\protect\restr{P^{-1}(U_{i})} \ar[r]^-{\protect\ptwo} 
\ar[dl]^-{k_{j}} & U_{ij}. \\
X\times U_{ij}\times X \ar@/_1.5pc/[urr]_{\operatorname{pr}_{2}}}
\end{equation*}
In particular, we claim that
\begin{equation}
  \label{eq:6}
  k_{j}\circ k_{i}^{-1}(x,u,y) = \bigl(\phi_{ij}(u)(x), u,
  \phi_{ij}(u)(y)\bigr). 
\end{equation}
To see this, consider
\begin{align*}
  k_{j}\circ k_{i}^{-1}(x,u,y) &=
  k_{j}\bigl(h_{i}^{-1}(u,x),h_{i}^{-1}(u,y)\bigr) \\
&= \bigl(\xi_{j}\bigl(h_{i}^{-1}(u,x)\bigr),u,
\xi_{j}\bigl(h_{i}^{-1}(u,y)\bigr)\bigr) 
\end{align*}
which equals the right-hand side of \eqref{eq:6} in view of
\eqref{eq:4}~and \eqref{eq:1}.  

We include the following lemma 
to motivate some of the constructions in the
next section.
At this point, it will be convenient to introduce the notation $\dw$
for $p(\omega)$.  This will make some of the more complicated formulas
below, and elsewhere, a bit easier to digest.

\begin{lemma}
  \label{lem-hm-1}
  Suppose that $p:\Omega\to Y$ is a locally trivial fibre bundle with
  fibre $X$ and that $G=\Omega\pstar\Omega$.  If $H$ is a topological
  groupoid and if $\alpha:G\to H$ is a
  continuous groupoid \hm, then there is open cover $\set{U_{i}}$ of
  $Y$ and continuous maps $\psi_{i}:p^{-1}(U_{i})\to H$ such that
  \begin{equation}
    \label{eq:8}
    \alpha(\omega,\omega')=\psi_{i}(w)\psi_{i}(\omega')^{-1}\quad
    \text{for $(\omega,\omega')\in G\restr{p^{-1}(U_{i})}$.}
  \end{equation}
\end{lemma}
\begin{remark}
  \label{rem-hm-2}
  It is implicit in the statement of the Lemma that
  $r\bigl(\psi_{i}(\omega
  )\bigr)=r\bigl(\alpha(\omega,\omega^{\prime})\bigr)$, while
  $s\bigl(\psi_{i}(\omega)\bigr)$ and, more generally,
  $\psi_{i}(\omega)^{-1}\psi_{j}(\omega)$ depend only on $p(\omega)$,
  for $\omega\in p^{-1}(U_{ij})$.
\end{remark}
\begin{proof}
  [Proof of Lemma~\ref{lem-hm-1}]
  Choose $\set{U_{i}}$, $h_{i}:p^{-1}(U_{i})\to U_{i}\times X$ and
  $k_{i} :G\restr{p^{-1}(U_{i})}\to X\times U_{i}\times X$ as above.
  For each $i$, fix $z_{i}\in X$ and define
  \begin{equation*}
    \psi_{i}(\omega):=
     \alpha\circ k_{i}^{-1}\bigl(\xi_{i}(\omega),\dw,z_{i}\bigr),
  \end{equation*}
and observe that if $(\omega,\omega')\in G\restr{p^{-1}(U_{i})}$, then
\begin{align*}
  \alpha(\omega,\omega') &:= \alpha\circ
  k_{i}^{-1}\bigl(\xi_{i}(\omega), \dw, \xi_{i}(\omega')\bigr) \\
&=\alpha\circ
  k_{i}^{-1}\bigl(\xi_{i}(\omega),\dw,z_{i}\bigr)\alpha\circ
  k_{i}^{-1}\bigl( z_{i},\dw,\xi_{i}(\omega')\bigr) \\
&= \psi_{i}(\omega)\psi_{i}(\omega')^{-1}
\end{align*}
as required.
\end{proof}

\section{The Dixmier-Douady Class}

First, it will be helpful to recast \cite{fmw:xx}*{Theorem~1} in terms
of Brauer groups as defined in \cite{kmrw:ajm98}.\footnote{We follow
  the notation and terminology of \cite{kmrw:ajm98}.  In particular,
  $\sbr(G)$ denotes the collection of pairs $({\A},\alpha)$ of the
  kind we have been considering, while $\br(G)$ denotes their Morita
  equivalence classes.}  To start off, we only need $G$ to be a proper
and not necessarily locally trivial.
\begin{thm} [cf.~\cite{kmrw:ajm98}*{Proposition~11.2}]
  \label{thm-iso}
  If $G$ is a second countable, locally compact, proper principal
  groupoid, then there is an isomorphism of $\brg$ onto $\br(\gomg)$
  which sends $(\A,\alpha)\in\sbrg$ to the class of the bundle
  representing the crossed product $\cs(G;\A)$ in $\sbr(\gomg)$.  The
  inverse is given by sending $\B\in\sbr(\gomg)$ to
  $(p^{*}\B,\tau\tensor1)$ where $p:\go\to\gomg$ is the quotient map,
  $p^{*}\B=\set{(\omega,b):\dw=p_{\B}(b)}$ is the pull-back and
  $(\tau\tensor 1)_{(\omega,\omega')}(\omega',b)=(\omega,b)$.
\end{thm}
\begin{proof}
  It follows from \cite{mrw:jot87}*{Proposition~2.2} and
  \cite{kmrw:ajm98}*{Theorem~4.1} that $[\A,\alpha]\mapsto [\A^{\go}]$
  is an isomorphism of $\brg$ onto $\br(\gomg)$ , where
  $\A^{\go}=s_{*}(\A)/G$ and where
  $s_{*}(\A)$ is the pull-back of $\A$ to $G$ via
  the source map $s$ (cf.~\cite[p.~914]{kmrw:ajm98}).  Furthermore,
  \cite{fmw:xx}*{Theorem~1} implies that $\bunsec(\gomg,\A^{\go})$ is
  Morita equivalent to $\cs(G;\A)$.  This proves the first assertion.
  
  On the other hand, \cite{kmrw:ajm98}*{Theorem~4.1} implies that the
  inverse is given by
  \begin{equation*}
    [\B]\mapsto [\B^{\goo},\id^{\goo}],
  \end{equation*}
where $\goo$ is $\go$ regarded as a left-$G$,
right-$\gomg$ equivalence (cf.~\cite[p.~924]{kmrw:ajm98}). 
  However, $\goo$ is easily seen to be isomorphic (as an equivalence)
  to the graph of $p$ $\gr(p)$, where $\gr(p):=\set{(\omega,\dw):\omega\in\go}$ (cf. \S
  6 of \cite{kmrw:ajm98}). Therefore \cite{kmrw:ajm98}*{Lemma~6.5}
  implies that
\begin{equation*}
  [\B^{\goo},\id^{\goo}]=[\B^{\gr(p)},\id^{\gr(p)}]=[p^{*}\B,p^{*}\id]. 
\end{equation*}
Noting that $p^{*}\id=\tau\tensor1$ completes the proof.
\end{proof}

As a corollary, we get the following version of
\cite{raeros:tams88}*{Theorem~1.1}.  
\begin{cor}
  \label{cor-rr-thm1.1}
  Suppose that $[\A,\alpha]\in\brg$.  Then there is a locally trivial
  $\K$-bundle $\bundlefont{B}$ over $Y$ such that
  \begin{equation*}
    [\A,\alpha]=[p^{*}\bundlefont{B},\beta].
  \end{equation*}
  In particular, the Dixmier-Douady class $\delta(\A)$ must lie in
  the image $p^{*}\bigl(H^{3}(Y;\Z)\bigr)$ of the map
  $p^{*}:H^{3}(Y;\Z)\to H^{3}(\Omega;\Z)$ induced by $p$.
\end{cor}
\begin{proof}
  Since $\tau_{\go}:=(\go\times \K,\tau\tensor1)$ acts as the identity
  in $\brg$, we can replace $(\A,\alpha)$ with
  $\bigl(\A\tensor_{\go}(\go\times\K),
  \alpha\tensor(\tau\tensor1)\bigr)$.  Thus we may as well assume that
  $A:=\bunsec(\go,\A)$ is stable.  It follows from
  Theorem~\ref{thm-iso} and \cite{raewil:tams85}*{Proposition~1.4(1)}
  that $\delta(\A):=\delta(A)\in p^{*}\bigl(H^{3}(Y;\Z)\bigr)$.  Thus
  there is a stable continuous-trace \cs-algebra $B$ with spectrum $Y$
  such that $A=p^{*}B$ \cite{raewil:tams85}*{Proposition~1.4(2)}.  In
  view of \cite{rw:morita}*{Proposition~5.59}, we may assume that
  $B=\bunsec(Y,\B)$ where $\B$ is a locally trivial $\K$-bundle with
  structure group $\Aut\K$ as required.
\end{proof}

Let $\B$ be a locally trivial $\K$-bundle over $Y$ as in
Corollary~\ref{cor-rr-thm1.1}.  It follows from
\cite{rw:morita}*{Proposition~4.53} for example, that there is a cover
$\U=\set{U_{i}}$ of $Y$ by \emph{pre-compact} sets and continuous
functions $\sigma_{ij}$ from $U_{ij}$ into $\Aut K$ with the
point-norm topology such that $\sigma=\set{\sigma_{ij}}$ form a
$1$-cocycle in the \emph{set} $Z^{1}(\U,\P)$ (where we'll use $\P$ for
the sheaf of germs of continuous $\Aut\K$ valued functions on $Y$) such that $\B$
is isomorphic to
\begin{equation}\label{eq:12}
  \coprod_{i} U_{i}\times \K/\!\!\sim,
\end{equation}
where we identify $(j,u,T)$ with $\bigl(i,u, \sigma_{ij}(u)(T)\bigr)$.
Furthermore, the Dixmier-Douady class $\delta(\B)$ is the image of
$\Delta([\sigma])\in H^{2}(Y,\sheaffont{S})$ in $H^{3}(Y;\Z)$ where
$\Delta:H^{1}(Y,\P)\to H^{2}(Y,\sheaffont{S})$ is the bijection given
in \cite{rw:morita}*{Proposition~4.83}.  Thus, refining $\U$ if
necessary, we can assume that there are continuous functions
$\mu_{ij}$ from $U_{ij}$ into the unitary group $\UK$ of $\K$ with the
strong topology such that
\begin{equation}
  \label{eq:11}
  \sigma_{ij}(u)=\Ad \mu_{ij}(u)\quad\text{for all $u\in U_{ij}$.}
\end{equation}
Then $\Delta([\sigma])$ is represented by $\nu:=\set{\nu_{ijk}}\in
H^{2}(\U,\sheaffont{S})$ given by
\begin{equation*}
  \nu_{ijk}(u)\mu_{ik}(u)=\mu_{ij}(u)\mu_{jk}(u).
\end{equation*}

If $\B$ is of the form \eqref{eq:12}, then $\A:=p^{*}\B$ is of the
form
\begin{equation}
  \label{eq:13}
  \coprod_{i}p^{-1}(U_{i})\times \K/\!\!\sim,
\end{equation}
where we identify $(j,\omega,T)$ with
$\bigl(i,\omega,\sigma_{ij}(\dw)(T)\bigr) $.  Thus if
$(p^{*}\B,\alpha)\in\sbrg$, there must be continuous groupoid \hm s
\begin{equation}
  \label{eq:14}
  \pi_{i}:G\restr{p^{-1}(U_{i})}\to\Aut\K
\end{equation}
such that
\begin{equation*}
  \alpha(\omega,\omega')\bigl([i,\omega',T]\bigr) = [i,\omega,
  \pi_{i}(\omega, \omega')(T)].  
\end{equation*}
Note that if $(\omega,\omega')\in G\restr{p^{-1}(U_{i})}$, then $\dw =
\dot \omega'$ and
\begin{equation}
  \label{eq:15}
  \sigma_{ij}(\dw)\circ \pi_{i}(\omega,\omega')=
  \pi_{j}(\omega,\omega') \circ \sigma_{ij}(\dw).
\end{equation}
Assume that $G$ is locally trivial with $k_{i}$, $h_{i}$, $\xi_{i}$
and $\phi_{ij}$ defined as in \S1.  Then we fix some $z_{i}\in X$ and
define
\begin{equation*}
  \psi_{i}:\piui\to\Aut\K
\end{equation*}
by
\begin{equation*}
  \psi_{i}(\omega)=\pi_{i}\bigl(k_{i}^{-1}\bigl(\xi_{i}(\omega),\dw,
  z_{i} \bigr)\bigr).
\end{equation*}
Then, as in Lemma~\ref{lem-hm-1},
\begin{equation*}
  \pi_{i}(\omega,\omega')=\psi_{i}(\omega)\psi_{i}(\omega')^{-1}.
\end{equation*}
Now we observe that
\begin{align}
  \label{eq:16}
  \psi_{i}(\omega)^{-1}\circ {}&\sigma_{ij}(\dw)\circ \psi_{j}(\omega)  \\
  &=\pi_{i}\bigl(k_{i}^{-1}\bigl(\xi_{i}(\omega),\dw,a_{i}\bigr)\bigr)^{-1}
  \circ
  \pi_{i}\bigl(k_{j}^{-1}\bigl(\xi_{j}(\omega),\dw,z_{j}\bigr)\bigr)
  \circ \sigma_{ij}(\dw) \notag \\
  \intertext{which, since
    $\phi_{ji}(\dw)(\xi_{j}(\omega))=\xi_{i}(\omega)$, is} &=
  \pi_{i}\bigl(k_{i}^{-1}\bigl(z_{i},\dw,\xi_{i}(\omega)\bigr)\bigr)
  \circ \pi_{i}\bigl(k_{i}^{-1}\bigl(\xi_{i}(\omega),\dw,
  \phi_{ji}(\dw)\bigr)\bigr) \circ \sigma_{ij}(\dw) \notag \\
  &= \pi_{i}\bigl(k_{i}^{-1}\bigl(z_{i},\dw,
  \phi_{ji}(\dw)(z_{j})\bigr)\bigr)\circ\sigma_{ij}(\dw).\label{eq:18}
\end{align}
In particular, \eqref{eq:16} depends only on the class $\dw$ of
$\omega$ in $Y$.
\begin{thm}
  \label{thm-calculation}
  Suppose that $G$ is a second countable, locally compact, locally
  trivial, proper principal groupoid, that $\B$ is a locally trivial
  $\K$-bundle over $Y$ and
  $(p^{*}\B,\alpha)\in\sbrg$, $\U$, $\psi_{i}$ and $\sigma_{ij}$ are
  as defined above.  Then
  \begin{equation}
    \label{eq:17}
    \beta_{ij}(\dw)=\psi_{i}(\omega)^{-1}\circ\sigma_{ij}(\dw)\circ
    \psi_{j}(\omega) 
  \end{equation}
  defines a class $\beta$ in $H^{1}(Y,\P)$ which depends only on
  $[p^{*}\B,\alpha] \in\brg$.  Furthermore, the Dixmier-Douady class
  of the corresponding groupoid crossed product $\cs(G;p^{*}\B)$ is
  the image of $\Delta(\beta)$ in $H^{3}(Y;\Z)$.
\end{thm}
\begin{proof}
  It follows from Theorem~\ref{thm-iso} that
  $\delta\bigl(\cs(G;p^*\B)\bigr)$ depends only on $[p^{*}\B,\alpha]$.
  It is clear from \eqref{eq:18} and \eqref{eq:17} that
  $\set{\beta_{ij}}$ defines a cocycle in $Z^{1}(\U,\P)$ and therefore
  a class $\beta$ in $H^{1}(Y,\P)$.  It will suffice to see that the
  image of $\Delta(\beta)$ is $\delta\bigl(\cs(G;p^{*}\B)\bigr)$.
  
  Since the proof of \cite{fmw:xx}*{Theorem~1} implies that
  $\cs(G,p^{*}\B)$ is Morita equivalent to $\bunsec(Y,p^{*}\B/G)$
  where $p^{*}\B/G$ is the orbit space on $p^{*}\B$ with respect to
  the right action of $G$ given by $[i,\omega,T]\cdot (\omega,\omega')
  = [i,\omega', \pi_{i}(\omega,\omega')(T)]$, it will suffice to
  compute $\delta(p^{*}\B/G)$.  We'll denote the image of
  $[i,\omega,T]$ in $p^{*}\B/G$ by $\pbgclass{i,\omega,T}$.  Let
  $\V=\set{V_{i}}$ be a cover of $Y$ such that $\overline V_{i}\subset
  U_{i}$ (\cite{rw:morita}*{Lemma~4.32}).  Let $P$ be a rank-one
  projection in $\K$.  If $w\in\piui$, then
  \begin{equation*}
    [i,\omega,\psi_{i}(\omega)(T)]\cdot(\omega,\omega') =
    [i,\omega',\pi_{i}(\omega,\omega') \circ\psi_{i}(\omega)(T)] =
    [i,\omega', \psi_{i}(\omega')(T)].
  \end{equation*}
  Therefore
\begin{equation*}
  \dw\mapsto \pbgclass{i,\omega,\psi_{i}(\omega)(P)}
\end{equation*}
is a rank-one projection field in $(p^{*}\B/G)\restr{\piui}$.  Fell's
vector-valued Tietze extension theorem
\cite{fd:representations}*{Theorem~II.14.8} implies that the
restriction to $\overline{V_{i}}$ extends to a global section
$q_{i}\in \bunsec(Y,p^{*}\B/G)$ such that
\begin{equation*}
  q_{i}(\dw)=\pbgclass{i,\omega,\psi_{i}(\omega)(P)}\quad\text{for all
  $\dw\in V_{i}$.}
\end{equation*}
We may as well assume that there are continuous functions
$\theta_{ij}:U_{ij}\to \UK$ such that
\begin{equation*}
  \beta_{ij}(\dw)=\Ad\theta_{ij}(\dw)\quad\text{for all $\dw\in U_{ij}$.}
\end{equation*}
Of course, $\V$ is a refinement of $\U$ with the same index set and
$\Delta(\beta)$ is given by the $2$-cocycle $\set{\epsilon_{ijk}}$ on
$\V$ defined by
\begin{equation*}
  \theta_{ij}(\dw)\theta_{jk}(\dw)=
  \epsilon_{ijk}(\dw)\theta_{ik}(\dw)\quad\text{for
  $\dw\in V_{ijk}$.} 
\end{equation*}
As above, there is a $v_{ij}\in\bunsec(Y,p^{*}\B/G)$ such that
\begin{equation*}
  v_{ij}(\dw)=\pbgclass{i,\omega,\psi_{i}(\omega)\bigl(P\theta_{ij}(\dw)^{*}
  \bigr)}
  \quad\text{for $\dw\in \overline{V_{ij}}$.}
\end{equation*}
For $\dw\in V_{ij}$ we certainly have
\begin{equation}
  \label{eq:19}
  v_{ij}(\dw)v_{ij}(\dw)^{*}=
  \pbgclass{i,\omega,\psi_{i}(\omega)(P)}=q_{i}(\dw)  ,
\end{equation}
while
\begin{equation}
  \label{eq:20}
  \begin{split}
    v_{ij}(\dw)^{*}v_{ij}(\dw) &=
    \pbgclass{i,\omega,\psi_{i}(\omega)\bigl(
      \theta_{ij}(\dw)P\theta_{ij}(\dw)^{*}\bigr)} \\
    &= \pbgclass{i,\omega,\psi_{i}(\omega)\circ\beta_{ij}(\dw)(P)} \\
    &= \pbgclass{i,\omega,\sigma_{ij}(\dw)\circ\psi_{j}(\omega)(P)} \\
    &= \pbgclass{j,\omega,\psi_{j}(\omega)(P)}\\
    &= q_{j}(\dw).
  \end{split}
\end{equation}
Thus we have just the set-up for \cite{rw:morita}*{Lemma~5.28}.  So to
find the Dixmier-Douady class, we compute that
\begin{align*}
  v_{jk}(\dw)^{*}v_{ij}(\dw)^{*} &=
  \pbgclass{j,\omega,\psi_{j}(\omega)\bigl( \theta_{jk}(\dw)P\bigr)}
  \pbgclass{i,\omega ,\psi_{i}(\omega)\bigl(\theta_{ij}(\dw)P\bigr)}
  \\
  &= \pbgclass{i,\omega,
    \sigma_{ij}(\dw)\circ\psi_{j}(\omega)\bigl(\theta_{jk}(\dw)P\bigr)}
  \pbgclass{i,\omega,\psi_{i}(\omega)\bigl(\theta_{ij}(\dw)P\bigr)} \\
  &=
  \pbgclass{i,\omega,\psi_{i}(\omega)\circ\beta_{ij}(\dw)\bigl(\theta_{jk}(\dw)
    P\bigr)}
  \pbgclass{i,\omega,\psi_{i}(\omega)\bigl(\theta_{ij}(\dw)P\bigr)} \\
  &=
  \pbgclass{i,\omega,\psi_{i}(\omega)\bigl(\theta_{ij}(\dw)\theta_{jk}(\dw)
    P\theta_{ij}(\dw)^{*}\bigr)}
  \pbgclass{i,\omega,\psi_{i}(\omega)\bigl(\theta_{ij}(\dw)P\bigr)} \\
  &= \epsilon_{ijk}(\dw)\pbgclass{i, \omega, \psi_{i}(\omega)\bigl(
    \theta_{ik}(\dw) P\bigr)} \\
  &= \epsilon_{ijk}(\dw) v_{ik}(\dw)^{*}.
\end{align*}
Thus \cite{rw:morita}*{Lemma~5.28} implies that
$\Delta(\beta)=[\set{\epsilon_{ijk}}]$ gives the Dixmier-Douady class
of $\cs(G;p^{*}\B)$ as claimed.
\end{proof}

\begin{definition}
  \label{def-loc-uni-pb}
  Let $G=\Omega\pstar\Omega$ be a locally trivial proper principal
  groupoid, and let $\B$ be a locally trivial $\K$-bundle with
  structure group $\Aut \K$.  Then an action of $\alpha$ of $G$ on
  $p^{*}\B$ is called a \emph{locally unitary pull-back} if there is a
  cover $\U$ of $Y$ and continuous functions
  \begin{equation*}
    v_{i}:\piui\to\UK\quad\text{such that}\quad\psi_{i}(\omega)=\Ad
    v_{i}(\omega) \quad\text{for all $\omega\in\piui$.}
  \end{equation*}
\end{definition}

\begin{remark}
  \label{rem-loc-uni-pb}
  Note that if $\alpha$ is as in Definition~\ref{def-loc-uni-pb}, then
  there are actually groupoid \hm s
  \begin{equation*}
    u_{i}:G\restr{\piui}\to\UK\quad\text{such that}\quad
    \pi_{i}(\omega,\omega') =\Ad u_{i}(\omega,\omega'),
  \end{equation*}
  where
  $u_{i}(\omega,\omega')=\psi_{i}(\omega)\psi_{i}(\omega')^{-1}$.
\end{remark}

Since
\begin{equation*}
  \sigma_{ij}(\dw)\circ \pi_{i}(\omega,\omega') =
  \pi_{j}(\omega,\omega') \circ \sigma_{ij}(\dw),
\end{equation*}
there are continuous \emph{groupoid \hm s}
\begin{equation}\label{eq:22}
  \lambda_{ij}:G\restr{\piui}\to\C\quad\text{such that}\quad
  \sigma_{ij}(\dw)\,\bar{ }\,\bigl(u_{j}(\omega,\omega')\bigr) =
  \lambda_{ij}(\omega,\omega') u_{j}(\omega,\omega')
\end{equation}
\begin{remark}
  \label{rem-pr-ob}
  We want to think of $\set{\lambda_{ij}}$ as a \emph{generalized
    Phillips-Raeburn obstruction} for $\alpha$.  It is possible to
  view \eqref{eq:22} as defining a class in a certain sheaf cohomology
  group, which we would naturally denote $H^{1}(Y,\sheaffont{G})$, but
  we will not pursue that here.
\end{remark}

% The class\altdb\footnote{We have to see the class is well-defined and
%   ``what it means''.} $\lambda$ in $H^{1}(Y,\sheaffont{G})$ is the
%   \emph{generalized Phillips-Raeburn obstruction} for $\alpha$.

  \begin{example}
    \label{ex-loc-uni-pr}
    Suppose that $\B$ is a locally trivial $\K$-bundle with transition
    functions $\set{\sigma_{ij}}$ as above and that $\gamma:H\to\Aut
    \bunsec(Y,\B)$ is a locally unitary automorphism group.  If
    $p:\Omega\to Y$ is a principal $H$-bundle, then let
    $G=\Omega\pstar\Omega$ and $\alpha$ the action of $G$ on $p^{*}\B$
    obtained by the pull-back: $\alpha:H\to\Aut p^{*}\B$ given by
    $\alpha_{s}(f)(\omega)= \gamma_{s}\bigl(f(s^{-1}\cdot
    \omega)\bigr)$. Since $\gamma$ is locally unitary, we can assume
    $\U$ has been taken so that there are strongly continuous \hm s
    $v^{i}:U_{i} \times H\to \UK$ such that
    \begin{equation*}
      \gamma(s\cdot\omega,\omega)\bigl([i,\dw,T]\bigr) = [i,\dw,
      v^{i}_{s}(\dw) T v^{i}_{s}(\dw)^{*}].
    \end{equation*}
    Then the continuous functions
    \begin{equation*}
      \zeta_{ij}:U_{ij}\to \widehat H
    \end{equation*}
    defined by
\begin{equation*}
  \sigma_{ij}(\dw)\bigl(v^{j}_{s}(\dw)\bigr) = \zeta_{ij}(\dw) (s)
  v^{i}_{s}(\dw)
\end{equation*}
determine a class $\zeta(\gamma)\in
H^{1}(Y,\widehat{\sheaffont{H}}_{\text{ab}})$ 
called the Phillips-Raeburn obstruction for~$\gamma$
(cf. \cite[Proposition~3.3]{echwil:jot01}).  (Here
$\widehat{\sheaffont{H}}_{\text{ab}}$ is the dual of the 
abelianization 
$H/\overline{[H,H]}$ of $H$ which coincides with the collection of
characters on $H$.)
Note that for
$\alpha$ as above,
\begin{equation*}
  \pi_{i}(\omega,\omega')=\Ad
  u_{i}(\omega,\omega')\quad\text{where}\quad u_{i}(s\cdot
  \omega,\omega) = v^{i}_{s}(\dw).
\end{equation*}
Thus the generalized Phillips-Raeburn obstruction and the original are,
in this case, related by
\begin{equation*}
  \lambda_{ij}(s\cdot\omega,\omega)=\zeta_{ij}(\dw)(s).
\end{equation*}
  \end{example}

  \begin{thm}
    \label{thm-main-loc-uni}
    Suppose that $G$ is a second countable, locally compact, locally
    trivial, proper principal groupoid, and that
    $(p^{*}\B,\alpha)\in\sbrg$ is a locally unitary pull-back, with
    $\B$ a locally trivial $\K$-bundle as above. Then, with
    $\lambda_{ij}$ as defined in \eqref{eq:22}, the equation
    \begin{equation*}
      \tau_{ijk}(\dw)=\lambda_{ij}\bigl(k_{j}^{-1}\bigl(z_{j},\dw,
      \phi_{kj}(\dw) (z_{k})\bigr)\bigr)
    \end{equation*}
    defines an element $\lip<\alpha,p>$ in $H^{3}(Y;\Z)$ which depends
    only on $[p^{*}\B,\alpha]\in\brg$ and $\delta(\B)$.  Furthermore,
    the Dixmier-Douady class of the crossed-product is given by the
    formula 
\begin{equation}\label{eq:21}
  \delta\bigl(\cs(G,p^{*}\B)\bigr)= \delta(\B)+\lip<\alpha,p>.
\end{equation}
  \end{thm}
  \begin{proof}
    Since the class of the crossed product only depends on
    $[p^{*}\B,\alpha]$, it will suffice to establish \eqref{eq:21}.
    We adopt the notation used in the proof of
    Theorem~\ref{thm-calculation}.  Thus
    $\delta\bigl(\cs(G;p^{*}\B)\bigr)$ is given by
    $\set{\epsilon_{ijk}}$ where
    \begin{gather*}
      v_{ik}(\dw)^{*}v_{ij}(\dw)^{*}=\epsilon_{ijk}(\dw)
      v_{ik}(\dw)^{*}\quad\text{for
        all $\dw\in V_{ijk}$, and} \\
      v_{ij}(\dw):=\pbgclass{i,\omega,\psi_{i}(\omega)(P)}.
    \end{gather*}
    We may assume that $\delta(\B)$ is given by a cocycle
    $\set{\nu_{ijk}}$ which is determined by functions
    $\mu_{ij}:U_{ij}\to\UK$ such that
\begin{equation*}
  \sigma_{ij}(\dw)=\Ad\mu_{ij}(\dw)\quad\text{and}\quad \nu_{ijk}(\dw)
  \mu_{ik}(\dw) = \mu_{ij}(\dw)\mu_{jk}(\dw).
\end{equation*}
Since $\alpha$ is a locally unitary pull-back, we can, in view of
\eqref{eq:18}, write
\begin{equation*}
  \beta_{ij}(\dw)=\Ad\theta_{ij}(\dw)\quad\text{where}\quad
  \theta_{ij}(\dw) = u_{i}\bigl(k_{i}^{-1}\bigl(z_{i},\dw,
  \phi_{ji}(\dw) (z_{j})\bigr)\bigr) \mu_{ij}(\dw)
\end{equation*}
with $u_{i}$ as defined in Remark~\ref{rem-loc-uni-pb}.  But then
\begin{align*}
  \theta_{ij}(\dw)\theta_{jk}(\dw) &=
  u_{i}\bigl(k_{i}^{-1}\bigl(z_{i},\dw,
  \phi_{ji}(\dw) (z_{j})\bigr)\bigr) \mu_{ij}(\dw) \\
  &\qquad\qquad
  u_{j}\bigl(k_{j}^{-1}\bigl(z_{j},\dw,\phi_{kj}(\dw)(z_{k})\bigr)\bigr)
  \mu_{jk} (\dw) \\
  &= \tau_{ijk}(\dw) u_{i}\bigl(k_{i}^{-1}\bigl(z_{i},\dw,
  \phi_{ji}(\dw) (z_{j})\bigr)\bigr) u_{i}
  \bigl(k_{j}^{-1}\bigl(z_{j},\dw,\phi_{kj}(\dw)(z_{k})\bigr)\bigr) \\
  &\qquad\qquad\qquad \mu_{ij}(\dw)\mu_{jk}(\dw) \\
  &= \tau_{ijk}(\dw) \nu_{ijk}(\dw) u_{i}\bigl(k_{i}^{-1}\bigl(
  z_{i},\dw,
  \phi_{ki}(\dw)(z_{k})\bigr)\bigr) \mu_{ik}(\dw). \\
  &= \tau_{ijk}(\dw) \nu_{ijk}(\dw)\theta_{ik}(\dw).
\end{align*}
It follows that $\set{\tau_{ijk}\nu_{ijk}}$ is a $2$-cocycle defining
$\delta\bigl(\cs(G,p^{*}\B)\bigr)$.  In particular, $\tau_{ijk}$ is
also a $2$-cocycle and we're done.
  \end{proof}

  \begin{example}
    \label{ex-last}
    We should compare our Theorem~\ref{thm-main-loc-uni}, in the
    setting of 
    Example~\ref{ex-loc-uni-pr}, with
    \cite{raeros:tams88}*{Theorem~1.5} .  To see that our formula
    coincides with Raeburn and Rosenberg's, we have to see that our
    $\lip<\alpha,p>$ is the same as their $\lip<\gamma,p>$.  However,
    \begin{equation*}
      \lambda_{ij}\bigl(k_{j}^{-1}\bigl(z_{j},\dw,\phi_{kj}(\dw)
      (z_{k}) \bigr)\bigr) = \lambda_{ij}\bigl(h_{j}^{-1}(z_{j},\dw),
      h_{j}^{-1}( z_{k}\cdot s_{kj}(\dw),\dw)\bigr).
    \end{equation*}
    If, as we may, we take $z_{i}=e$ for all $i$, then
\begin{equation*}
  \lambda_{ij}\bigl(\omega,s_{kj}(\dw)\cdot\omega
  \bigr)=\zeta_{ij}(\dw)(s_{jk}(\dw)) 
\end{equation*}
which certainly appears to be a representative for Raeburn and
Rosenberg's $\lip<\gamma, p>$.  Unfortunately we used a different
convention in \eqref{eq:1} and \eqref{eq:10} from the one used in
\cite{raeros:tams88}*{Theorem~1.5}.  In view of this, our
$\lip<\alpha,p>$ should be $-\,\lip<\gamma,p>$.  The difference seems to
be that the expression ``$\Ad\bigl(u^{i}_{\lambda_{ij}(t)}(t)
v_{ij}(t)\bigr)$'' in the middle of page~15 
in the proof of 
\cite{raeros:tams88}*{Theorem~1.5} should be
``$\Ad\bigl(u^{i}_{\lambda_{ji}(t)}(t) v_{ij}(t)\bigr)$''.
  \end{example}

\bibliographystyle{amsxport} 
\bibliography{fmw-ddc}

\end{document}